\documentclass[twocolumn, a4paper]{article}

\usepackage[T1]{fontenc}
\usepackage[latin9]{inputenc}
\usepackage{amsthm}
\usepackage{amsmath}
\usepackage{amssymb}
\usepackage{graphicx}
\usepackage{cite}

\newcommand{\revision}[1]{#1}

\title{Reduction of interaction delays in networks\\
	{\small[appeared in EPL, 103 (2013) 10006, DOI:10.1209/0295-5075/103/10006]}}

\author{L. L\"ucken$^1$, J.P. Pade$^1$, K. Knauer$^2$ and S. Yanchuk$^1$}

\date{}

\begin{document}

\maketitle

\emph{$^1$ Humboldt-University of Berlin, Institute of Mathematics, Unter den Linden 6, 10099 Berlin, Germany.}\\
\emph{$^2$ Universit{\'e} Montpellier 2, Institut de Math{\'e}matiques et de Mod{\'e}lisation de Montpellier (I3M), Place Eug{\`e}ne Bataillon, 34095 Montpellier Cedex, France.}

\begin{abstract}
Delayed interactions are a common property of coupled natural systems
\revision{and therefore arise in a variety of different applications}.
For instance, signals in neural or laser networks propagate at finite
speed giving rise to delayed connections. Such systems are often modeled
by delay differential equations with discrete delays. 
\revision{In realistic situations, these delays are not identical on 
different connections.} We show that by a componentwise timeshift 
transformation it is often possible to
reduce the number of different delays and simplify the
models \revision{without loss of information}. We identify dynamic invariants of this transformation, determine
its capabilities to reduce the number of delays and interpret these
findings in terms of the topology of the underlying graph. 
In particular, we show that networks with identical 
\revision{sums of delay times} along
the fundamental semicycles are dynamically equivalent and \revision{we} provide
a normal form for these systems. We illustrate the theory using a
network motif of coupled Mackey-Glass systems with 8 different time
delays, which can be reduced to an equivalent motif with three delays.
\end{abstract}

\section{Introduction}

Time delays play an important role in \revision{the realistic modeling of} 
networks of dynamical systems.
For instance, in gene regulatory networks delayed couplings
arise due to the finite duration of biochemical reaction chains~\cite{Danino2010,Sevim2010,OBrien2012,Li2007,Heiden1979,Mackey1977a},
\revision{in population dynamics due to the processes of maturation 
and gestation~\cite{Kuang1993},}
and in laser networks \revision{they correspond to} the propagation time of light between
interacting lasers~\cite{Soriano2013,DHuys2008,Franz2008,Nixon2012}.
Similarly, in neuronal networks
delays occur due to the finite propagation times of action potentials
along the axons or to reaction times at chemical synapses\revision{~\cite{Leibold2005,Rossoni2005,Memmesheimer2006,Ko2007,Bakkum2008,Campbell2007}}. Depending
on the physiological properties of axons and synapses, the delays
of signals between different cells in the network differ and in practice
there are as many different delays in the network as there are links
between single neurons. Such a diversity of delays was shown to enable
the design of robust pattern generating systems~\cite{Yanchuk2011,Popovych2011}
and to play an important role in neural processing of spatio-temporal
information~\cite{Carr1993}. \revision{On a mesoscopic level, 
networks of neuronal populations responsible for motor control and human 
cognition involve transmission delays depending on the locations of the 
populations in the brain~\cite{Noori2010,Lisman2005}.} Hence, to study such natural systems,
one would need to study models involving numerous different delays.
Independently of the physical origin, dynamics in networks of $N$
interacting systems can often be described by the equations
\begin{equation}
\dot{\boldsymbol{x}}_{j}(t)=f_{j}\left(\boldsymbol{x}_{j}(t),\boldsymbol{x}_{1}\left(t-\tau_{j1}\right),\dots,\boldsymbol{x}_{N}\left(t-\tau_{jN}\right)\right),\label{eq:first}
\end{equation}
where $\boldsymbol{x}_{j}(t)\in\mathbb{R}^{n_{j}}$, $1\le j\le N$,
denotes the state of the node $j$ and $\tau_{jk}\ge0$ is the time
delay, which occurs due to the signal propagation time from a node
$k$ to the node $j$. For instance, $\boldsymbol{x}_{j}$ may be
considered as a variable describing the dynamics of one neuron in
the network and $\tau_{jk}$ as the time which an action potential
needs to travel along the axon of neuron $\boldsymbol{x}_{k}$ until
it arrives at the synapse connecting it to $\boldsymbol{x}_{j}$.
Equation~(\ref{eq:first}) includes the possibility of an all-to-all
coupling structure. \revision{However,} in many situations, the coupling is not global
and $\boldsymbol{x}_{j}$ receives direct inputs from only relatively
few other nodes. Therefore, we consider the set $P_{j}$ of all 
nodes of the network, which are connected to the node $j$
and write system (\ref{eq:first}) in a more compact, equivalent form
\begin{equation}
\dot{\boldsymbol{x}}_{j}(t)=f_{j}\left(\boldsymbol{x}_{j}(t),\left(\boldsymbol{x}_{k}\left(t-\tau_{jk}\right)\right)_{k\in P_{j}}\right).\label{eq:one}
\end{equation}
\revision{The only difference of the form of writing (\ref{eq:one}) from 
the general form (\ref{eq:first}) is that is explicitly shows
which connections are present in the network.}
Throughout what follows, we consider directed networks, thus one may
generally have $\tau_{jk}\ne\tau_{kj}$ or even $k\in P_{j}$ while
$j\notin P_{k}$. Moreover, we assume that the network is connected
and hence, the number $L$ of links is at least $N-1$. In fact, the
obtained results hold also for networks with multiple links between
the nodes, but we discuss 
\revision{only systems of the form (\ref{eq:one})} in order to avoid
complicated notations. For a network of $N$ systems, up to $N^{2}$
different time delays may occur in (\ref{eq:one}). 
This creates immense challenges for the analysis of the
system, since every single delay may alter the properties and dynamics
significantly~\revision{\cite{Erneux2009,Foss1996,Giacomelli1996,Hale2000,Flunkert2010,Heiligenthal2011,Kanter2011}}.

In this Letter we show how the number of different delays can be reduced.
It appears that any connected network (\ref{eq:one}) possesses a
characteristic number of delays which are essential for describing
the dynamics. This number of essential delays equals the cycle space
dimension $C=L-\left(N-1\right)$ of the underlying graph and is usually
smaller than the number of distinct $\tau_{jk}$. We show that the
essential delays correspond to \revision{the sums of the delay times} along fundamental
semicycles (cycles in the underlying undirected graph) in the network.
\revision{As we explain below,}
networks which have the same local dynamics and the same set of essential
delays can be considered as equivalent from the dynamical point of
view.
\revision{Our} results represent fundamental theoretical insights which 
help to understand delay coupled systems.
For instance, they can be used to simplify bifurcation analysis when varying one or 
several delays \revision{or} to speed up numerical simulations of delay equations. 

\subsection{\revision{Reduction of delays in unidirectional rings}}
As a simple illustrative example, let us consider a ring of unidirectionally
coupled systems 
\begin{equation}
\dot{\boldsymbol{x}}_{j}\left(t\right)=f\left(\boldsymbol{x}_{j}\left(t\right),\boldsymbol{x}_{j+1}\left(t-\tau_{j}\right)\right),\label{eq:ring-inhom}
\end{equation}
with inhomogeneous coupling delays $\tau_{j}$, $\ 1\le j\le N$,
and periodic indices, i.e., $\boldsymbol{x}_{N+1}=\boldsymbol{x}_{1}$.
It is known~\cite{Sande2008,Perlikowski2010} that this system can
be reduced to a ring where all time delays are equal to $\tau=\frac{1}{N}\sum\tau_{j}$.
This can be done by a componentwise timeshift transformation for
each node 
\begin{equation}
\boldsymbol{y}_{j}\left(t\right)=\boldsymbol{x}_{j}\left(t+\eta_{j}\right),\quad1\le j\le N.\label{eq:timeshift-transformation}
\end{equation}
Indeed, the transformed variables $\boldsymbol{y}_{j}\left(t\right)$
fulfill 
\begin{equation}
\dot{\boldsymbol{y}}_{j}\left(t\right)=f\left(\boldsymbol{y}_{j}\left(t\right),\boldsymbol{y}_{j+1}\left(t-\tau_{j}+\eta_{j}-\eta_{j+1}\right)\right).\label{eq:transformed-ring}
\end{equation}
\revision{This means that the timeshift transformation (\ref{eq:timeshift-transformation}) formally allows to
change the delays as $\tau_{j}\mapsto\tilde{\tau}_{j}=\tau_{j}-\eta_{j}+\eta_{j+1}$.}
If the shifts $\eta_{j}$ are chosen appropriately, all hitherto
distinct delays become the same while their sum along the ring, the
roundtrip time $N\tau$, is preserved.
It is clear that both systems, (\ref{eq:ring-inhom}) and (\ref{eq:transformed-ring}),
possess the same steady states and there is a one-to-one correspondence
of periodic solutions \revision{via (\ref{eq:timeshift-transformation})}. 
However, \revision{previous studies using this type of transformation} did not attempt
to show in which sense (\ref{eq:one}) and (\ref{eq:transformed-ring})
are dynamical equivalent nor did they explore how the method generalizes
to more complex networks. Our aim is to treat these questions and
apply the timeshift transformation (\ref{eq:timeshift-transformation})
to reduce the number of different delays in networks with arbitrary
topology.

\section{Reduction of delays in general coupling structures}
Applying the transformation (\ref{eq:timeshift-transformation}) to \revision{a system of the form}
(\ref{eq:one}), we obtain the new system 
\begin{equation}
\dot{\boldsymbol{y}}_{j}\left(t\right)=f_{j}\left(\boldsymbol{y}_{j}\left(t\right),\left(\boldsymbol{y}_{k}\left(t-\tilde{\tau}_{jk}\right)\right)_{k\in P_{j}}\right)\label{eq:transformed-system}
\end{equation}
with \revision{the modified delays} $\tilde{\tau}_{jk}=\tau_{jk}+\eta_{k}-\eta_{j}$. As for the
unidirectional ring (\ref{eq:ring-inhom}),
also in arbitrary networks the roundtrip time\revision{s}
along closed paths remain unchanged by the componentwise timeshift
transformation. In order to formulate this general principle more
precisely, we \revision{use} two definitions. A \emph{semicycle} in a directed
graph is a \revision{set} $c=\revision{\{}\ell_{1},...,\ell_{k}\revision{\}}$ of links
$\ell_{j}$, which constitute a cycle when neglecting their directionality.
\revision{For example, in the network displayed in fig.~\ref{fig:basal-ganglia-nw}
the links along the path STN$\rightarrow$GPe$\rightarrow$STN
constitute a semicycle, but also those of GPe$\rightarrow$STN$\rightarrow$GPi$\leftarrow$GPe.}
The \emph{roundtrip time} along \revision{a} semicycle $c$ is
defined as 
\begin{equation}
T\left(c\right):=\Big|\sum_{j=1}^{k}\sigma_{j}\tau\left(\ell_{j}\right)\Big|,\label{eq:def-roundtrip}
\end{equation}
where $\tau\left(\ell_{j}\right)$ is the time delay along the link
$\ell_{j}$ and $\sigma_{j}$ is either $+1$ or $-1$ depending on
the direction of the link $\ell_{j}$ with respect to an arbitrary
but fixed orientation of the semicycle. 
\revision{For example, the roundtrip of the semicycle
GPe$\overset{\tau_2}{\rightarrow}$STN$\overset{\tau_3}{\rightarrow}$GPi$\overset{\tau_4}{\leftarrow}$GPe
in fig.~\ref{fig:basal-ganglia-nw} is $T=|\tau_2 + \tau_3 - \tau_4|$.}
Note that the orientation
can be chosen in two different ways, which corresponds to opposite
signs of $\sigma_{j}$. However, the obtained value for $T\left(c\right)$
is independent of the choice of orientation. Most importantly, the
roundtrip time stays invariant under the componentwise timeshift transformation,
i.e., 
\begin{equation}
\tilde{T}\left(c\right)=\Big|\sum_{j=1}^{k}\sigma_{j}\tilde{\tau}\left(\ell_{j}\right)\Big|=T\left(c\right),
\end{equation}
where $\tilde{\tau}(\ell)$ is the transformed delay on the link $\ell$.
Denoting the source of $\ell$ by $s(\ell)$ and its target by $t(\ell)$,
we have $\tilde{\tau}(\ell)=\tau(\ell)+\eta_{s(\ell)}-\eta_{t(\ell)}$. \revision{On the other hand, 
for any choice of delays $\tilde{\tau}(\ell)$ which preserves the roundtrips,
there exist corresponding timeshifts $\eta_j$ such that $\tilde{\tau}(\ell)=\tau(\ell)+\eta_{s(\ell)}-\eta_{t(\ell)}$.
This equivalence
provides an intuitive way of looking for possible transformations in simple networks:
One can redistribute delays as long as the roundtrips stay fixed.}

\subsection{Construction of a delay-free spanning tree}

In the following, we explain how to construct timeshifts $\eta_{j}$
such that the transformed network (\ref{eq:transformed-system}) has a 
minimal number of delayed connections. Simultaneously with the timeshifts, we 
construct a set $S=\left\{ \ell_{1},...,\ell_{N-1}\right\} $
of $N-1$ links, such that the transformed delays $\tilde{\tau}(\ell)$
are zero on all links from $S$ and non-negative on all other links.
In fact, the set $S$ is a spanning tree: a set of $N-1$ links which
does not contain any semicycle. 
Each node of the network occurs at least once as a target
or source of some link in $S$, and the addition of any link $\ell\notin S$
results in a set $S^{\prime}=S\cup\left\{ \ell\right\}$, which contains
exactly one semicycle. This is the fundamental semicycle corresponding
to $\ell$ \revision{with respect to $S$}. The transformed delay $\tilde{\tau}\left(\ell\right)$
on the link $\ell$ equals the roundtrip time of this fundamental
semicycle. Consequently, the number of delays in the transformed system
can be reduced to
\begin{equation}
C=L-N+1,\label{eq:cycle-space-dim}
\end{equation}
where $L$ is the number of links in the network. The resulting system, whose links
are instantaneous along $S$, may function as a normal form for the
class of equivalent systems which exhibit the same roundtrips along
their semicycles.

\begin{figure}
\centering{}\includegraphics[width=0.75\linewidth]{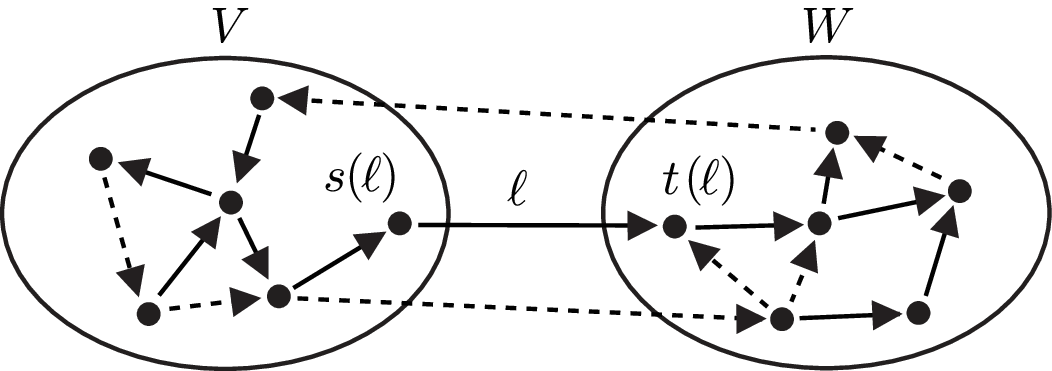}\caption{\label{fig:fundamental-section}Illustration of the partition of the
network into two sets $V$ and $W$, which is induced by a link $\ell$
contained in a spanning tree $S$ (solid links). Links which are not
contained in the spanning tree are indicated by dashed lines.}
\end{figure}

To construct the appropriate timeshifts and \revision{the spanning tree} $S$, we proceed iteratively.
Firstly, we select a spanning tree $S_{0}=\{\ell_{1}^{0},...,\ell_{N-1}^{0}\}$
in the following way: $\ell_{1}^{0}$ is chosen as a link with the
minimal delay and successive elements $\ell_{j}^{0}$ are always chosen
to have minimal possible delay, such that $\{\ell_{1}^{0},...,\ell_{j}^{0}\}$
remains cycle free. As a next step consider the link $\ell\in S_{0}$
with minimal non-zero delay $\tau(\ell)>0$ (it equals $\ell_{1}^{0}$
if $\tau(\ell_{1}^{0})>0$). As depicted in fig.~\ref{fig:fundamental-section},
the link $\ell$ divides $S_{0}$ in two connected parts: one either being empty 
or containing at least one link connecting to the source $s(\ell)$ of $\ell$,
the other being empty or connecting its target $t(\ell)$ with other nodes.
Correspondingly, the nodes are also divided into two sets $V$ and $W$ 
with $s(\ell)\in V$ and $t(\ell)\in W$.
 An application of the timeshift transformation with $\eta_{j}=\tau(\ell)$
for all $j\in W$, and $\eta_{j}=0$ for $j\in V$, yields $\tilde{\tau}(\ell)=\tau(\ell)+\eta_{s(\ell)}-\eta_{t(\ell)}=0$.
All other delays on links in $S_{0}$ remain unchanged
\revision{because for links which connect nodes from the same set $V$ or $W$ there is no difference in
the timeshifts and therefore no change of delays.}
By construction of $S_{0}$\revision{, any link $\ell^{\prime}\notin S_0$, which connects $V$ and $W$,
must satisfy $\tau(\ell^{\prime}) \ge \tau(\ell)$. Otherwise, $\ell^{\prime}$
had been added to $S_0$ instead of $\ell$. Therefore, we have
$\tilde{\tau}(\ell^{\prime}) = \tau(\ell^{\prime})+\eta_{s(\ell^{\prime})}-\eta_{t(\ell^{\prime})} \ge \tau(\ell^{\prime}) - \tau(\ell)\ge 0$. Hence, delays on links which are not contained in $S_{0}$ cannot become negative.}
A successive spanning tree $S_{1}$,
which contains at least one more instantaneous link is obtained by
keeping all instantaneous links from $S_{0}$ and then adding links
in the same manner as above. Proceeding in this fashion produces a
completely instantaneous spanning tree $S$ in less than $N$ steps.

\begin{figure}
\centering{}\includegraphics[width=0.6\linewidth]{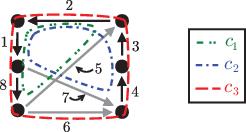}\caption{\label{fig:fundamental-cycle} A network of $N=6$ systems with cycle
space dimension $C=3$. A spanning tree {[}black links{]},
three fundamental links {[}grey links{]} and three corresponding
fundamental cycles {[}see legend{]}, are shown. The numbers refer
to the indexing of the links which is used in the text. According
to our results, any admissible dynamics in this network involves at
most three different interaction delays.}
\end{figure}

\revision{After the above reduction, the number of different delays
is not larger than $C=L-N+1$. Since $C$ is the minimal number of
different delays which can be achieved generically,}
we call \revision{it} the \emph{essential number of delays}.
\revision{Here, "generically" means that the set of delays, which allow
to reduce the number of different delays to less than $C$, is a set of 
measure zero in the parameter space $\mathbb{R}^L_{\ge0}$.}
$C$ has a well known meaning for network graphs~\cite{Diestel2010}:
It is the cycle space dimension which is defined as the maximal number
of independent cycles. A set of cycles is called independent if none
of the cycles can be obtained as a \revision{sum} of the 
other \revision{cycles} in the set. \revision{Here, the sum of two cycles 
$c_1$ and $c_2$ is defined as the symmetric difference 
$c_1 \bigtriangleup c_2 =(c_1 \cup c_2) \setminus (c_1 \cap c_2)$.}
For instance, in the network given
in fig.~\ref{fig:fundamental-cycle} the cycle space dimension is $C=3$.

\revision{
\subsection{A two-cycle neuronal network}
Let us illustrate the main ideas of the reduction with an examplary
network structure which resembles the connections between several 
areas in the brain (see fig.~\ref{fig:basal-ganglia-nw}):
the subthalamic nucleus (STN), the external segment (GPe), and 
the internal segment (GPi) of the globus pallidus \cite{Terman2002, Tachibana2011}.
Additional external inputs from the thalamus and cortex (ThCx) arrive at the STN.
We denote the different areas by $\boldsymbol{x}_{1}$ (STN), $\boldsymbol{x}_{2}$ (GPe),
$\boldsymbol{x}_{3}$ (GPi), and $\boldsymbol{x}_{4}$ (ThCx).
As general form of the equations which may describe the corresponding dynamics,
we assume
\begin{equation}
\begin{aligned}
\dot{\boldsymbol{x}}_{1}\left(t\right) & =  f_{1}\left(\boldsymbol{x}_{1}\left(t\right),\boldsymbol{x}_{2}\left(t-\tau_{2}\right),\boldsymbol{x}_{4}\left(t-\tau_{5}\right)\right),\\
\dot{\boldsymbol{x}}_{2}\left(t\right) & =  f_{2}\left(\boldsymbol{x}_{2}\left(t\right),\boldsymbol{x}_{1}\left(t-\tau_{1}\right)\right),\\
\dot{\boldsymbol{x}}_{3}\left(t\right) & =  f_{3}\left(\boldsymbol{x}_{3}\left(t\right),\boldsymbol{x}_{1}\left(t-\tau_{3}\right),\boldsymbol{x}_{2}\left(t-\tau_{4}\right)\right),\\ 
\dot{\boldsymbol{x}}_{4}\left(t\right) & =  f_{4}\left(\boldsymbol{x}_{4}\left(t\right),\boldsymbol{x}_{4}\left(t-\tau_6\right)\right),
\end{aligned}\label{eq:bg-nw}
\end{equation}
with connection delays $\tau_{j}>0$.
\begin{figure}
\revision{
\centering{}\includegraphics[width=0.85\linewidth]{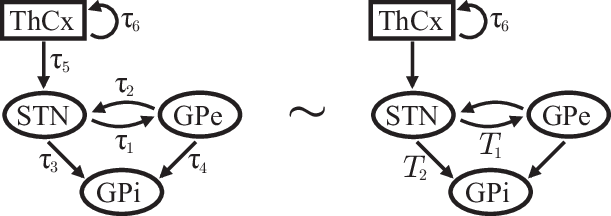} \caption{\label{fig:basal-ganglia-nw}
A small network with four nodes and two
cycles corresponding to eqs.~(\ref{eq:bg-nw}) which resemble the connectivity
between several regions in the brain [see main text]. The componentwise timeshift transformation (\ref{eq:timeshift-transformation}) 
leads from the original system (left) to a system with a reduced number of delays (right), where
$T_{1}=\tau_{1}+\tau_{2}$ and $T_{2}=\tau_{2}+\tau_{3}-\tau_{4}$.}}
\end{figure}

The transformation (\ref{eq:timeshift-transformation}) leads to new
delays 
\begin{multline*}
\tilde{\tau}_{1}=\tau_{1}-\eta_{2}+\eta_{1},\ \tilde{\tau}_{2}=\tau_{2}-\eta_{1}+\eta_{2}, \ \tilde{\tau}_{3}=\tau_{3}-\eta_{3}+\eta_{1},\\
 \tilde{\tau}_{4}=\tau_{4}-\eta_{3}+\eta_{2},\ \tilde{\tau}_{5}=\tau_{5}-\eta_{1}+\eta_{4},\text{ and }\tilde{\tau}_{6}=\tau_{6}.
\end{multline*}
The roundtrips $T_{1}:=\tau_{1}+\tau_{2}$ and $T_{2}:=\tau_{2}+\tau_{3}-\tau_{4}$
are invariant: they stay the same for the new delays.
In particular, this implies that a reduction to one delay, which was possible for the ring, 
cannot be achieved in general. There exist 6 formal possibilities 
to reduce the four different delays 
$\tau_{1-4}$ of (\ref{eq:bg-nw}) to two by an appropriate 
choice of $\eta_{1-3}$ [see Table~\ref{tab:Formal-reductions-for-motif2}].
However, some reductions should be excluded because they lead to negative delays in the resulting
network. Such ''anticipating arguments'' create functional differential equations
of mixed type. We avoid this case because for these equations
the initial value problem is usually ill-posed~\cite{Hale1993}.
As a result of the reduction, an equivalent system [fig.~\ref{fig:basal-ganglia-nw}, right panel]
contains only three time delays: the self-feedback $\tau_{6}$  to the ThCx, the short roundtrip 
delay $T_{1}$ between STN and GPe, and the roundtrip $T_{2}$ along the STN $\to$ GPe $\to$ GPi connection.
The example illustrates also several general principles. Firstly the delays of self-couplings 
cannot be eliminated (e.g. $\tilde{\tau}_{6}=\tau_{6}$). Further, one may always choose 
$\eta_{4}=\eta_{1}-\tau_{5}$, such that $\tilde{\tau}_{5}=0$. That is, delays on leaves 
(nodes with only one link to others) may be neglected.
}

\begin{table}
\revision{
\begin{centering}
\begin{tabular}{|c|c|c|c|c|c|}
\hline 
$\tilde{\tau}_{1}$  & $\tilde{\tau}_{2}$  & $\tilde{\tau}_{3}$  & $\tilde{\tau}_{4}$ & $\tilde{\tau}_{5}$ & $\tilde{\tau}_{6}$\tabularnewline
\hline 
\hline 
$T_{1}$  & 0  & 0  & $-T_{2}$ & 0 & $\tau_6$\tabularnewline
\hline 
$T_{1}$  & 0  & $T_{2}$  & 0 & 0 & $\tau_6$\tabularnewline
\hline 
0  & $T_{1}$  & $T_{2}-T_{1}$  & 0 & 0 & $\tau_6$\tabularnewline
\hline 
0  & $T_{1}$  & 0  & $T_{1}-T_{2}$ & 0 & $\tau_6$\tabularnewline
\hline 
$T_{1}/2$  & $T_{1}/2$  & $T_{2}-T_{1}/2$  & 0 & 0 & $\tau_6$\tabularnewline
\hline 
$T_{1}/2$  & $T_{1}/2$  & 0  & $T_{1}/2-T_{2}$ & 0 & $\tau_6$\tabularnewline
\hline 
\end{tabular}
\par\end{centering}

\caption{\label{tab:Formal-reductions-for-motif2}Possible delay reductions
in the network {[}fig.~\ref{fig:basal-ganglia-nw}, eqs.~(\ref{eq:bg-nw}){]}
under the condition $\tilde{\tau}_5=0$. The two roundtrips
$T_{1}=\tau_{1}+\tau_{2}$ and $T_{2}=\tau_{2}+\tau_{3}-\tau_{4}$
were used to parametrize the different combinations for $\tilde{\tau}_{1-4}$.
The spanning tree algorithm {[}see main text{]} produces one of the first 
four possibilities depending on the values of $\tau_{1-4}$.}
}
\end{table}

\revision{Let us give a brief explanation of how the algorithm works
in the present example with delay times $\tau_1=1$, $\tau_2=2$, $\tau_3=3$, $\tau_4=4$, and $\tau_5=5$
For this choice we find that $T_1=3$ and $T_2=1$.
Let us denote the links by $\ell_{1-5}$ such that $\tau(\ell_j)=\tau_j$.
The initial spanning tree is then $S_0=\{\ell_1,\ell_3,\ell_5\}$. Note that,
although $\tau_2<\tau_3$, the link $\ell_2$ cannot be added to $S_0$ instead 
of $\ell_3$ because this would result in a cycle $\{\ell_1,\ell_2\}$. Similarly,
$\ell_5$ is added instead of $\ell_4$.
In the next step, the sets $V = \{\boldsymbol{x}_1,\boldsymbol{x}_3,\boldsymbol{x}_4\}$ and $W=\{\boldsymbol{x}_2\}$
are formed and the delay $\tau_1$ is eliminated by a timeshift transformation
with $\eta_1=\eta_3=\eta_4=0$ and $\eta_2=\tau_1=1$. This results in transformed delays
$\tilde{\tau}_1=0$, $\tilde{\tau}_2=\tilde{\tau}_3=3$, and $\tilde{\tau}_4=\tilde{\tau}_5=5$.
The successive spanning tree has to be chosen as before $S_1=S_0$. The smallest non-zero
delay on $S_1$ is now $\tilde{\tau}_3=3$ and correspondingly we set
$V = \{\boldsymbol{x}_1,\boldsymbol{x}_2,\boldsymbol{x}_4\}$ and $W=\{\boldsymbol{x}_3\}$.
An application of the timeshift transformation with $\eta_1=\eta_2=\eta_4=0$ and $\eta_3=3$
eliminates the delay on $\ell_3$ and we obtain:
$\tilde{\tilde{\tau}}_1=\tilde{\tilde{\tau}}_3=0$, $\tilde{\tilde{\tau}}_2=3$,
$\tilde{\tilde{\tau}}_4=2$, and $\tilde{\tilde{\tau}}_5=5$.
For the final step, we find $S_2=S_0$, $V=\{\boldsymbol{x}_4\}$, $W=\{\boldsymbol{x}_{1-3}\}$
and, applying the timeshift-transform with $\eta_4=0$ and $\eta_{1-3}=5$ gives
the reduced delay distribution which is listed in the fourth row of
Table~\ref{tab:Formal-reductions-for-motif2}.}

\section{Dynamical invariants}
As mentioned above, it is possible to create negative delays $\tilde{\tau}_{jk}$
by certain choices of the timeshifts $\eta_{j}$. This would change
the type of the problem fundamentally. At least for this reason it
is worth considering to what extent the systems (\ref{eq:one}) and
(\ref{eq:transformed-system}) possibly differ in their dynamical
behavior. Fortunately, negative delays turn out to be the only pathological
case, as we show in the following. More precisely, we show that characteristic
exponents of steady states and periodic solutions are the same in
both systems. This implies that the linear stability of these solutions
is the same~\cite{Hale1993}. Firstly, we consider a steady state
$\boldsymbol{x}^{\ast}=(\boldsymbol{x}_{1}^{\ast},\dots,\boldsymbol{x}_{N}^{\ast})$
of (\ref{eq:one}). It is as well a steady state of the transformed
system (\ref{eq:transformed-system}). For each characteristic exponent
$\lambda$ of $\boldsymbol{x}^{\ast}$ in (\ref{eq:one}) there is
a solution $\boldsymbol{\chi}(t)=\exp(\lambda t)\boldsymbol{\chi}_{0}$
of the linearized equation 
\begin{equation}
\dot{\boldsymbol{\chi}}_{j}(t)=\partial_{1}f_{j}^{\ast}\boldsymbol{\chi}_{j}(t)+\sum_{k\in P_{j}}\partial_{k+1}f_{j}^{\ast}\boldsymbol{\chi}_{k}(t-\tau_{jk}),\label{eq:vareq-1}
\end{equation}
where we introduce $\partial_{1}f_{j}^{\ast}:=\partial_{\boldsymbol{x}_{j}(t)}f(\boldsymbol{x}^{\ast})$
and $\partial_{k+1}f_{j}^{\ast}:=\partial_{\boldsymbol{x}_{k}(t-\tau_{jk})}f(\boldsymbol{x}^{\ast})$,
$1\le j\le N$~\cite{Hale1993}. A straightforward calculation shows
that the function $\exp(\lambda t)\tilde{\boldsymbol{\chi}}_{0}$
with $\tilde{\boldsymbol{\chi}}_{0,j}=\exp(\lambda\eta_{j})\boldsymbol{\chi}_{0,j}$
is a solution of the linearized equation for $\boldsymbol{x}^{\ast}$
in (\ref{eq:transformed-system}). Thus, $\lambda$ is a characteristic
exponent for the corresponding steady state in (\ref{eq:transformed-system})
as well. Further, consider a periodic solution $\boldsymbol{x}^{\ast}(t)=\boldsymbol{x}^{\ast}(t+p)$
of (\ref{eq:one}) with period $p>0$. According to the Floquet theory
for systems with time delays \cite{Hale1993}, a characteristic (Lyapunov)
exponent $\lambda$ of $\boldsymbol{x}^{\ast}(t)$ corresponds to
a solution $\boldsymbol{\chi}(t)=\exp(\lambda t)\boldsymbol{q}(t)$
of the linearized equation (\ref{eq:vareq-1}) with time-periodic
coefficients $\partial_{1}f_{j}^{\ast}(t)=\partial_{\boldsymbol{x}_{j}(t)}f(\boldsymbol{x}^{\ast}(t))$
and $\partial_{k+1}f_{j}^{\ast}(t)=\partial_{\boldsymbol{x}_{k}(t-\tau_{jk})}f(\boldsymbol{x}^{\ast}(t))$,
$1\le j\le N$, and a $p$-periodic function $\boldsymbol{q}(t)=(q_{1}(t),...,q_{N}(t))$.
A corresponding solution $\tilde{\boldsymbol{\chi}}(t)=\exp(\lambda t)\tilde{\boldsymbol{q}}(t)$
of the variational equation of the transformed solution in (8) is
obtained with $\tilde{q}_{j}(t)=\exp(\lambda\eta_{j})q_{j}(t+\eta_{j})$.

Using theoretical tools from the theory of semidynamical systems,
it is possible to show that a strong equivalence between (\ref{eq:one})
and (\ref{eq:transformed-system}) holds~\cite{Luecken2013a}. 
\revision{In particular, there is a natural one-to-one correspondence between invariant sets 
of the original system (\ref{eq:one}) and of the transformed
system (\ref{eq:transformed-system}). Any invariant set of (\ref{eq:transformed-system}) 
consists of the timeshifted trajectories of a corresponding set of (\ref{eq:one}) and \emph{vice versa}.
Moreover, maximal Lyapunov exponents of the corresponding sets coincide and they have
the same type of stability. The more abstract semidynamical systems approach 
can extend the results to systems with more 
general types of local dynamics as, for example, delay coupled partial
differential equations.}

\subsection{Motif of coupled Mackey-Glass systems}

The following example is intended to provide an idea of possible applications
and merits of the results presented in this article. Let us consider
a network motif of coupled Mackey-Glass systems~\cite{Mackey1977a,DHuys2008,Appeltant2011,Soriano2012}
\begin{equation}
\dot{x}_{j}\left(t\right)=-\gamma x_{j}\left(t\right)+\beta\frac{c\sum x_{k}(t-\tau_{jk})}{1+\left(c\sum x_{k}(t-\tau_{jk})\right)^{10}},\label{eq:MackeyGlassNW}
\end{equation}
where the summation is performed over all coupled nodes $k\in P_{j}$.
The coupling topology is shown in fig.~\ref{fig:fundamental-cycle}.
The parameters are fixed to $\gamma=0.1$, $c=0.525$, and $\beta=0.2$.
Although the total number of different delays is 8, only three roundtrips
$T\left(c_{j}\right)$, $j=1,2,3$, along the fundamental cycles indicated
in fig.~\ref{fig:fundamental-cycle} are important. Therefore, the
system can be reduced to an equivalent system with only three essential
delays $T(c_{1})=\tau_{5}+\tau_{2}+\tau_{1}+\tau_{8},$ $T(c_{2})=\tau_{7}+\tau_{4}+\tau_{3}+\tau_{2}+\tau_{1}$,
and $T(c_{3})=\tau_{6}+\tau_{4}+\tau_{3}+\tau_{2}+\tau_{1}+\tau_{8}$. 

\begin{figure}
\centering{} \includegraphics[width=0.9\columnwidth]{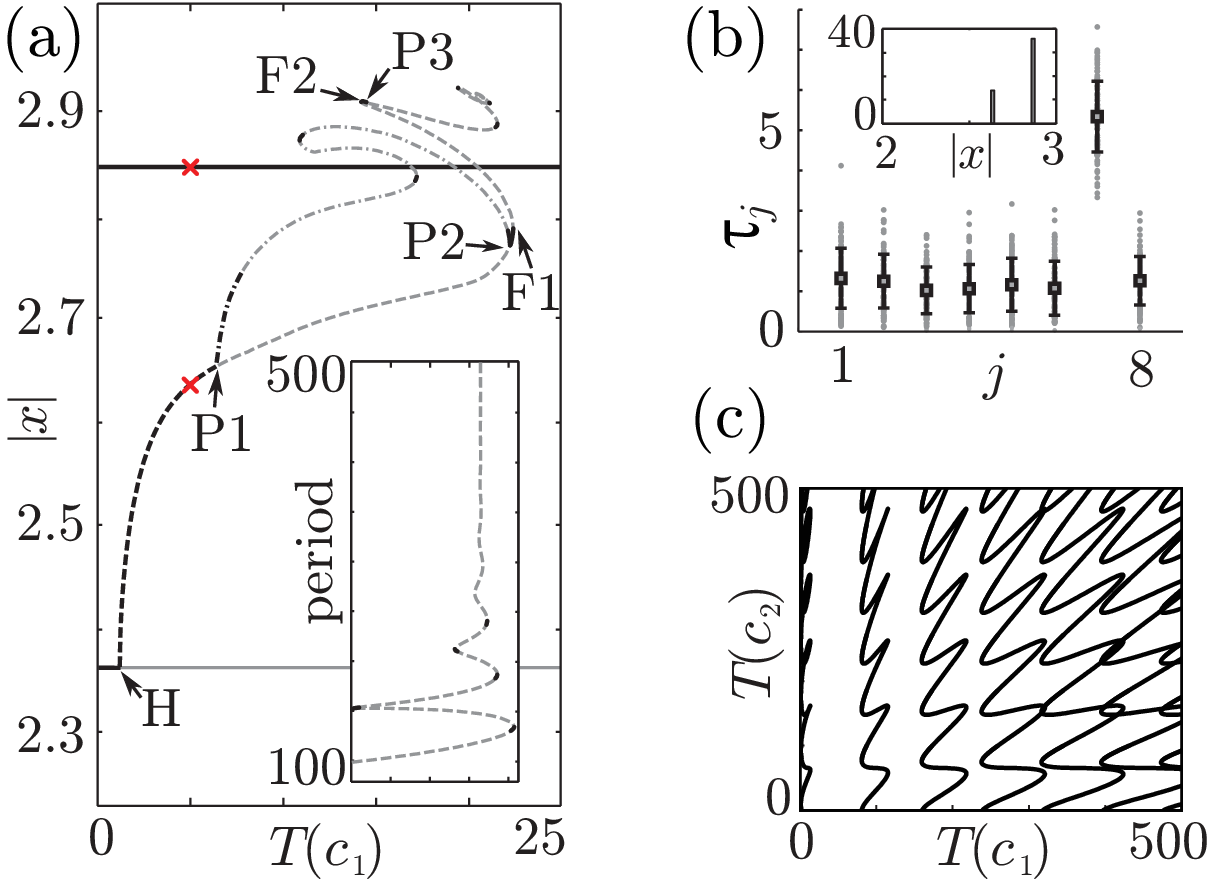}
\caption{\label{fig:Dynamical-Illustration}Dynamical features of six-node
network motif from fig.~\ref{fig:fundamental-cycle} of Mackey-Glass
systems (\ref{eq:MackeyGlassNW}). (a): Bifurcation diagram for $T\left(c_{1}\right)\in\left[0,23\right]$
and fixed $T\left(c_{2}\right)=10$, $T\left(c_{3}\right)=7$. Black
lines correspond to stable solutions and gray to unstable. Solid branches
are equilibria, dashed are periodic solutions. Bifurcations are denoted
by 'H' (Hopf), 'P{[}$n${]}' (period doubling), and 'F{[}$n${]}'
(fold). The ordinate measures the norm $|x|$ and the inset shows
the period along the main branch. Red crosses indicate the observed
solutions from chart (b). (b): 50 random realizations of single delays
which keep constant $T\left(c_{1}\right)=5$, $T\left(c_{2}\right)=10$,
and $T\left(c_{3}\right)=7$. Each dot corresponds to a realization
of $\tau_{j}$. Mean and standard deviation are indicated by squares
and bars. For each choice of delays, constant initial functions were
drawn from a uniform distribution on $\left[0,1.5\right]$. The inset
bars indicate the number of observed solutions $\boldsymbol{x}\left(t\right)$,
$t\in[10^{5},1.1\times10^{5}]$, with the norm $|\boldsymbol{x}|$
in their base interval. Chart (c) shows the trace and reappearances
of the Hopf bifurcations in the $T\left(c_{1}\right)$-$T\left(c_{2}\right)$-plane. }
\end{figure}

Firstly, the delays in the network were set up randomly in such a
way that the roundtrips were preserved with $T\left(c_{1}\right)=5$,
$T\left(c_{2}\right)=10$, and $T\left(c_{3}\right)=7$ {[}fig.~\ref{fig:Dynamical-Illustration}(b){]}.
For each realization of delays, the system was integrated numerically
with some randomly chosen initial conditions {[}see caption to fig.~\ref{fig:Dynamical-Illustration}{]}.
As a result of each integration, an attractor was found: either a
steady state or a periodic solution. The steady state was the same
for all choices of individual delays, while the periodic solutions
were connected by timeshift transformations of the form (\ref{eq:timeshift-transformation}).
All corresponding periodic solutions exhibit the same componentwise
supremum norm $\left|\boldsymbol{x}\right|=\sqrt{\sum_{j=1}^{N}\sup_{t}\left|x_{j}\left(t\right)\right|^{2}}$ {[}see
inset of fig.~\ref{fig:Dynamical-Illustration}(b){]}, which is invariant
under the transformation (\ref{eq:timeshift-transformation}). Thus,
we observe that the dynamical features of the system do not change
provided the roundtrips $T(c_{j})$ stay fixed, even though we have
chosen the individual delays randomly for each integration.

Now let us assume that one essential delay $T(c_{1})$ is a parameter
in the system, and two others stay fixed as $T\left(c_{2}\right)=10$
and $T\left(c_{3}\right)=7$. The bifurcation diagram with respect
to $T\left(c_{1}\right)$ {[}fig.~\ref{fig:Dynamical-Illustration}(a){]}
was obtained using numerical continuation software \cite{Engelborghs2002}.
We observe that for $T\left(c_{1}\right)=0$ two stable equilibria
are attractors for the system. The lower equilibrium looses stability
in a supercritical Hopf bifurcation (H) and the emergent periodic
solution undergoes a period doubling bifurcation (P1) which connects
to a reverse period doubling at (P2). The winding of the main branch
accompanied by period doubling bridges and divergent periods indicates
the presence of a homoclinic saddle-focus~\cite{Glendinning1984}.
A step further can be done by a two-parametric bifurcation diagram
{[}fig.~\ref{fig:Dynamical-Illustration}(c){]}, where two round-trips
$T(c_{1})$ and $T(c_{2})$ are changed. There, the Hopf bifurcation
(H) was followed numerically in the region $0\le T\left(c_{1}\right)\le500$
and $0\le T\left(c_{2}\right)\le500$. 

The role of the obtained bifurcation diagrams is that they describe
the dynamics of the system using the proper time delay parameters.
For instance, any combinations of individual delays $\tau_{1},\dots,\tau_{8}$,
which preserve the same roundtrip $T(c_{3})=7$, would represent a
point in a two dimensional bifurcation diagram in fig.~\ref{fig:Dynamical-Illustration}(c).
Therefore, the knowledge of the proper parameters allowed to reduce
the number of effective parameters from 8 to 3. In particular, such
a reduction of the dimensionality of the parameter space makes the
thorough study of the possible dynamical regimes in the considered
system feasible. 

\section{Summary and discussion}

We have presented a universal tool for the reduction of discrete
interaction-delays in networks of dynamical systems. It allows
to simplify the coupling structure without any loss of information. We have shown
that the number of essential delays which determine the dynamics
equals the cycle space dimension of the network, and that systems 
whose roundtrip times along the fundamental semicycles coincide are 
dynamically equivalent. 
\revision{Even if the roundtrip times are given, the distribution of the delays which
yields them is not unique. This provides the opportunity to choose
a delay distribution that is best suitable for the  problem under consideration.}
In some cases, as for instance for a unidirectional ring
of identical systems with inhomogeneous delays, 
\revision{the right choice of timeshifts} allows 
to uncover hidden symmetries \revision{and thereby to simplify the 
analysis of the system considerably~\cite{Baldi1994,Perlikowski2010}.
Conversely, the same transformation was applied 
to design pattern generators~\cite{Popovych2011,Yanchuk2011}
starting from the system with homogeneous delays.}
\revision{Our results are} an important step
towards the systematization and classification of dynamics on networks
with multiple delays. As a main consequence for general coupling structures, 
the bifurcation analysis can be simplified to a great extent as we have shown 
in the example of coupled Mackey-Glass systems. The benefits of the reduction seem to be the highest for networks 
with sparse connectivity. Nevertheless, in highly connected networks, 
reductions of delays along selected links (e.g. the ones with highest 
connection weight) or in separate motifs can be a useful application. 
Furthermore, it seems a promising endeavor to study the potential of 
the delay transformation to speed up numerical procedures. We have 
observed that the time it takes to integrate numerically two equivalent 
systems may vary extremely depending on the distribution of the individual 
delays. Finally, the results which were presented in this paper can be 
extended straightforwardly to many other types of local dynamics like 
time-discrete systems or evolution equations.

\subsubsection*{Acknowledgments}
We thank M. Zaks for useful discussions and the
DFG for financial support in the framework of the Collaborative Research
Center (SFB) 910 and the International Research Trainig Group (IRTG) 1740.

\bibliographystyle{unsrt}
\bibliography{Ybibliography.bib}

\end{document}